\documentclass[12pt]{iopart}

\usepackage{iopams}

\expandafter\let\csname equation*\endcsname\relax

\expandafter\let\csname endequation*\endcsname\relax

\usepackage{amsmath,bm,graphicx}

\begin{document}

\title[]{Stiff-response-induced instability for chemotactic bacteria and flux-limited Keller-Segel equation}

\author{Ben\^oit PERTHAME$^1$ and Shugo YASUDA$^2$}

\address{$^1$Sorbonne Universit\'es, UPMC Univ Paris 06, Laboratoire Jacques-Louis Lions  UMR CNRS 7598,  Universit\'e Paris-Diderot, Inria de Paris, F75005 Paris}
\address{$^2$Graduate School of Simulation Studies, University of Hyogo, Kobe 650-0047, Japan}
\ead{$^1$benoit.perthame@upmc.fr}
\ead{$^2$yasuda@sim.u-hyogo.ac.jp}

\vspace{10pt}
\begin{indented}
\item[]
\end{indented}

\begin{abstract}
	Collective motion of chemotactic bacteria as {\it E. Coli} relies, at the individual level, on a continuous reorientation by runs and tumbles. 
		It has been established that the length of run is decided by a stiff response to a temporal sensing of chemical cues along the pathway. 
		
		We describe a novel mechanism for pattern formation stemming from the stiffness of chemotactic response relying on a kinetic chemotaxis model which includes a recently discovered formalism for the bacterial chemotaxis. 
		We prove instability both for a microscopic description in the space-velocity space and for the macroscopic equation, a flux-limited Keller-Segel equation, which has attracted much attention recently.
		
		A remarkable property is that the unstable frequencies remain bounded, as it is the case in Turing instability. 
		Numerical illustrations based on a powerful Monte Carlo method show that the stationary homogeneous state of population density is destabilized and periodic patterns are generated in realistic ranges of parameters. 
		These theoretical developments are in accordance with several biological observations.
\end{abstract}

%
%
%
%
%

\section{Introduction}
Collective motion of chemotactic bacteria as {\it E. Coli} relies, at the individual level, on a continuous reorientation by runs and tumbles sensing extracellular chemoattractants produced by themselves \cite{art:75A,art:91BB,art:95BB,book:03B}. 
It has been established that the length of run is decided by a stiff response to temporal sensing of chemical cues along pathway, i.e., bacteria reduce their tumbling frequency and extend the run length as they sense an increase in concentrations of chemoattractants along the pathway.
Thus, the modulation of tumbling frequency in the chemotactic response is an essential mechanism for bacterial communities self-organization. 

This paper is concerned with the pattern formation of the population density of run-and-tumble chemotactic bacteria as {\it E. Coli}.
We describe a novel self-organized pattern formation mechanism stemming from a modulation of tumbling frequency with stiffness in chemotactic response.
Our analysis relies on a solid mathematical analysis and simulations using a unique Monte Carlo code.

In order to investigate the multiscale nature in this new self-organized pattern formation mechanism, we rely on a mesoscopic description, i.e., a kinetic reaction-transport equation for the chemotactic bacteria coupled with a reaction-diffusion equation for the chemoattractants.
The microscopic dynamic properties such as tumbling rate, modulation in stiff response, and proliferation (division/death) rate are included at the individual level.
We consider the following three main ingredients in the pattern formation: 
(i) the random run-and-tumble motion of bacteria, where the bacteria run linearly with a constant speed when rotating their flagella in counter-clockwise direction, but occasionally change the running directions (tumbling) when rotating their flagella in clockwise direction; 
(ii) the stiff and bounded signal response to the logarithmic sensing of chemoattractants along the pathway of bacterium, which generates the biased random motion searching for the higher-concentration region of chemoattractant; 
and (iii) the division/death of bacteria, where the population-growth rate depends on the local population density of bacteria. 
The kinetic transport equation considered in this paper describes all these ingredients at the microscopic (individual) level.

The pattern formations in the chemotaxis with population growth have been investigated at the macroscopic level by the Keller-Segel type equations \cite{art:71KS,art:71KS2,art:08TMPA}; for example, in \cite{art:91MMWM,art:96MT,art:07NPR,art:11PH}, the pattern formations induced by the properties of chemotaxis, i.e., so-called {\it chemotaxis-induced instability}, are demonstrated both theoretically and numerically.
Our paper is also concerned with the chemotaxis-induced instability, but which is based on the kinetic transport equation, which up to our knowledge, has not been carried out so far.

The kinetic approach has a distinctive advantage in studying the multiscale mechanism and mathematical hierarchy between the individual dynamics and macroscopic phenomena.
It has a long history and was first proposed in \cite{art:80A,art:88ODA} and then further developed toward involving the spatiotemporal variation of the chemoattractant along the pathway of bacterium \cite{art:05DS}, the internal dynamics of the cellular states \cite{art:04EO}, and the multi-cellular interactions \cite{art:07BBNS,art:08BLM}.
The mathematical foundations for the kinetic chemotaxis model have been strengthened involving the mathematical hierarchy between kinetic and continuum equations and the existence of solution for the kinetic chemotaxis equation \cite{art:00HO,art:02OH,art:04CMPS,art:08BCGP,art:08BC,art:10V,art:13JV,art:14STY,art:15AEV,art:15X,art:15L,art:16PTV,art:16C,art:17MW}.
The numerical methods to solve the kinetic chemotaxis equations have been also developed in \cite{art:10V,art:13RS,art:13RS2,art:14YF,art:16RKS,art:16BJP,art:17Y}.
The use of the kinetic chemotaxis model is also advanced due to the development of experimental technologies, which allow experimentalists to measure the individual velocities and turning angles in the collective motions of bacteria and give access to time scale measurements.
For example, in \cite{art:11SCBPBS,art:13RFWBYAMB,art:16EGBAV}, the advantage of the kinetic modeling is demonstrated by the comparison of the numerical and experimental results.

Our analysis also applies to the flux-limited Keller-Segel system, which is a very active research subject nowadays\cite{art:13JV,art:09PD,art:12CKWW,art:15BBTW,art:17BW}.
The flux-limited Keller-Segel system is derived as the asymptotic limit of the kinetic chemotaxis model mentioned earlier in the so-called ``diffusion limit'' \cite{art:PVW}. 
It incorporates a saturation of the chemotactic flux which avoids the blow-up of solutions. 
When diffusion is ignored, it has the property of finite speed of propagation. 
Thus the flux-limited Keller-Segel system can describe collective behaviours observed in various biological systems more realistically. 
Our instability result is also a new observation for the flux-limited Keller-Segel system.

In this paper, we propose a new mechanism leading to the linear instability of a kinetic chemotaxis equation coupled with a diffusion-reaction equation for the chemoattractant.
The kinetic chemotaxis equation involves a population-growth term, which depends on the local population density of bacteria, as well as a chemotactic response function in the integral kernel, which depends on the spatiotemporal variation of the chemoattractant along the pathway of each bacterium. 
We obtain a linear instability condition based on the stiffness of the response.
The stationary homogeneous state of the population density of bacteria becomes linearly unstable and periodic patterns are generated. 
We also numerically demonstrate the pattern formations by Monte Carlo method in which we vary the parameters involved in the linear instability condition.
A theoretical foundation for the Monte Carlo method is also presented.

\section{Main result}
Since the observation of the run-and-tumble movement of bacteria \cite{art:72BB,art:95B}, the kinetic chemotaxis equation has been proposed as an accurate description \cite{art:88ODA,art:04EO,art:05DS}.
In this study, we include a recently advocated formalism for bacterial chemotaxis, i.e., a logarithmic sensing \cite{art:09KJTW} and a stiff and bounded signal response \cite{art:83BSB}.
That is,
\begin{equation}\label{eq_kinetic}
	\begin{split}
		&\partial_tf(t,\bm{x},\bm{v})+\bm{v}\cdot \nabla f=\\
		&\frac{1}{k}\left\{\frac{1}{4\pi}\int_VK[D_t \log S|_{\bm{v}'}]f(\bm{v}')d\Omega(\bm{v}')
		-K[D_t\log S|_{\bm{v}}]f(\bm{v})\right\}
		+P[\rho]f(\bm{v}).
	\end{split}
\end{equation}
Here $f(t,\bm{x},\bm{v})$ is the microscopic population density of bacteria with a velocity $\bm{v}\in V$ at position $\bm{x}\in \mathbb{R}$ and time $t\ge 0$.
The velocity space of $\bm{v}$, $V$ is the surface of the unit ball, i.e., $|\bm{v}|=1$ and $\Omega(\bm{v})$ is the unit measure on $V$.
The tumbling kernel $K[D_t\log S|_{\bm{v}}]$ represents the stiff and bounded signal response (which is explained in Eq. (\ref{eq_kernel})) to the logarithmic sensing of chemical attractant $S(t,{\bm x})$ along the pathway with velocity ${\bm v}$.
Here $D_t X|_{\bm{v}}$ is the material derivative, i.e., $D_t X|_{\bm{v}}=\partial_t X+\bm{v}\cdot\nabla X$.  
The concentration of the chemical attractant $S(t,\bm{x})$ and macroscopic population density of bacteria $\rho(t,\bm{x})$ are calculated as, respectively,
\begin{equation}\label{eq_S}  
	-d\Delta S(t,\bm{x})+S(t,\bm{x})=\rho(t,\bm{x}),
\end{equation}
\begin{equation}\label{eq_rho}
	\rho(t,\bm{x})=\frac{1}{4\pi}\int_V f(t,\bm{x},\bm{v})d\Omega(\bm{v}),
\end{equation}
where $d$ is the molecular diffusion constant.
In Eq.~(\ref{eq_kinetic}), $P[\rho]$ is the population-growth rate of bacteria which depends on the local population density $\rho$ as
\begin{subequations}\label{eq_P}  
	\begin{align}
		&P[\rho]>0,\quad \mbox{for } 0<\rho<1,\\
		&P[\rho]<0,\quad \mbox{for } \rho>1,\\
		&P[\rho]\simeq 1-\rho,\quad \mbox{for }\rho\simeq 1.
	\end{align}
\end{subequations} 
Thus, the bacteria may divide when the local population density is lower than unity and the new born bacteria have the same velocities as the parents, but they may die when the local population density is larger than the unity.

We remark that the existence of the traveling wave in the kinetic transport equations with population growth (but without chemotactic responses, i.e., $K[X]=$const.) are proved in \cite{art:00H,art:00S,art:15BCN}, where new born particles may choose their velocities according to a prescribed equilibrium velocity distribution. 
In this paper, we use the simplest population-growth model among those for which the existence of the traveling wave is proved mathematically and the logistic population-growth term is recovered in the continuum limit.

The tumbling kernel $K[X]$ in Eq.~(\ref{eq_kinetic}) is a decreasing function of $X$ and we choose it as 
\begin{equation}\label{eq_kernel}
	K[X]=1-F[X],
\end{equation}
where $F[X]$ is a smooth and bounded function which satisfies the following properties,
\begin{subequations}\label{cond_F}
	\begin{align}  
		&F[0]=0,\\
		&\frac{d F[X]}{dX}>0,\\
		&F[X]\rightarrow \pm \chi \mbox{ as } X\rightarrow \pm \infty.
	\end{align}
\end{subequations}
Here, $F[X]$ represents the chemotactic response of the bacteria, say the response function, and $\chi$ represents the amplitude of modulation in the chemotactic response and takes a constant value between $0\le\chi<1$.

In Eqs.~(\ref{eq_kinetic})--(\ref{eq_P}), all quantities are nondimensionalized by the following characteristic quantities; 
the characteristic time $t_0$ is defined as $t_0=|(\frac{d\tilde P}{d\rho})_{\rho=1}|^{-1}$, where $\tilde P$ represents the population-growth rate in the dimensional form, the characteristic length $L_0$ is defined as $L_0=t_0V_0$, where $V_0$ is a constant speed of the bacteria.
The nondimensional parameter $k$ is defined as $k=1/(t_0\psi_0$), where $\psi_0$ is a mean tumbling frequency of the bacteria. The population density is scaled by that in the uniform stationary state $\rho_0$ and the concentration of the chemoattractant is scaled by $(a_0/b_0)\rho_0$, where $a_0$ is the production rate of chemoattractant by the bacteria and $b_0$ is the degradation rate of chemoattractant. 

It is easily seen that Eqs.~(\ref{eq_kinetic})--(\ref{eq_P}) have a constant uniform solution with $f(t,\bm{x},\bm{v})=1$, $S(t,\bm{x})=1$, and $\rho(t,\bm{x})=1$.
In our main result, this uniform solution gives to Turing-like instability \cite{art:52T,book:15P}. That is, the uniform solution is linearly unstable if the stiffness of the response function $F'[0]$ is sufficiently large as
\begin{equation}\label{kinetic_instability}
	\frac{F'[0]}{k}>\inf_{\lambda}\left[
		1+\frac{k}{\frac{k\lambda}{\arctan(k\lambda)}-1}
	\right]
	(1+d\lambda^2).
\end{equation}
In addition, the unstable eigenmodes are bounded, i.e., no high frequency oscillations occur and thus patterns are formed.

Furthermore, Eq.~(\ref{kinetic_instability}) includes the linear stability condition of a flux-limited Keller-Segel equation obtained by the asymptotic analysis of the kinetic chemotaxis equation Eq.~(\ref{eq_kinetic}) under a diffusion scaling introduced in Eq.~(\ref{eq_scale}).
This proves that our instability condition is sharp in the continuum limit ($k\rightarrow 0$).

\begin{figure}[p]
	\centering
	\includegraphics*[scale=0.8]{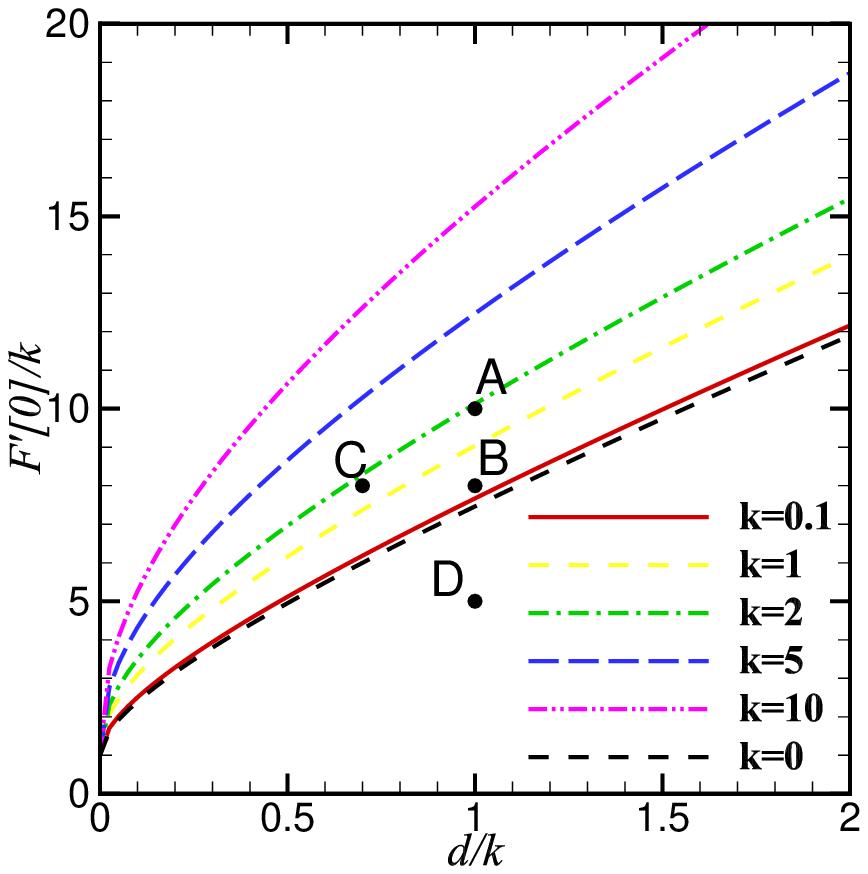}
	\caption{The diagram of the kinetic instability obtained by Eq.~(\ref{kinetic_instability}). The linear instability takes place when the stiffness of the response function $F'[0]/k$ exceeds the critical line of each $k$. 
	The numerical simulations are performed for the marks A, B, C, and D with several values of $k$. See Table \ref{t_para}.
}\label{fig_dki}
\end{figure}

Figure \ref{fig_dki} shows the linear instability diagram.
We numerically calculate the minimum values of the right-hand side of Eq.~(\ref{kinetic_instability}) with variation in the Fourier mode $\lambda$ for given values of $k$ and $d$.
One can observe that as increasing the stiffness of the response function $F'[0]$, the stationary homogeneous state with $f=S=1$ is destabilized (the {\it chemotaxis-induced} instability).
The critical lines increase monotonically as the diffusion coefficient $d$ increases, so that the stationary homogeneous state is more destabilized when the diffusion coefficient $d$ becomes smaller.
It is also seen that as decreasing $k$, the critical line decreases but converges to that for the flux-limited Keller-Segel equation.
The stationary homogeneous state is more likely destabilized as decreasing $k$.

\section{Linear instability analysis}
\subsection{Linearization}
It is easily seen that Eq.~(\ref{eq_kinetic}) has the uniform solution, $f(t,\bm{x},\bm{v})=S(t,\bm{x})=1$.
We carry out the linear instability analysis about this uniform solution.
We consider a small perturbation around the uniform solution as
\begin{subequations}\label{eq_perturbation}
	\begin{align} 
		f(t,\bm{x},\bm{v})&=1+ g(\bm{x},\bm{v})e^{\mu t},\\
		S(t,\bm{x})&=1+ S_g(\bm{x})e^{\mu t},\\
		\rho(t,\bm{x})&=1+ \rho_g(\bm{x})e^{\mu t},
	\end{align}
\end{subequations}
where $\mu$ is a constant which distinguishes stability ($\mathrm{Re}(\mu)<0$) and instability ($\mathrm{Re}(\mu)>0$).
From Eqs.~(\ref{eq_P})--(\ref{cond_F}), we can linearize the population-growth rate $P[\rho]$ and tumbling kernel $K[D_t\log S|_{\bm{v}}]$ in Eq. (\ref{eq_kinetic}) as, respectively,
\begin{equation}
	P[1+\rho_ge^{\mu t}]=-\rho_g e^{\mu t},
\end{equation}
and
\begin{equation}
	K[D_t\log(1+S_g e^{\mu t})|_{\bm{v}}]=
	1-F'[0](\mu S_g+\bm{v}\cdot\nabla S_g)e^{\mu t}.
\end{equation}
Thus, Eq.~(\ref{eq_kinetic}) is linearized as
\begin{subequations}\label{eq_linear}
	\begin{align}
		k(\mu g(\bm{x},\bm{v})&+\bm{v}\cdot\nabla g(\bm{x},\bm{v}))e^{\mu t}\nonumber \\
		&=\frac{1}{4\pi}\int_V(1-F'[0](\mu S_g(\bm{x})+\bm{v}'\cdot\nabla S_g(\bm{x})) e^{\mu t}+g(\bm{x},\bm{v}')e^{\mu t})d\Omega(\bm{v}')\nonumber\\
		&-(1-F'[0](\mu S_g(\bm{x})+\bm{v}\cdot\nabla S_g(\bm{x})) e^{\mu t}+g(\bm{x},\bm{v})e^{\mu t})
		-k\rho_g(\bm{x}),\\
		k\mu g(\bm{x},\bm{v})&+k\bm{v}\cdot\nabla g(\bm{x},\bm{v})
		\nonumber\\
		&=\frac{1}{4\pi}\int_V g(\bm{x},\bm{v}')d\Omega(\bm{v}')
		+F'[0]\bm{v}\cdot\nabla S_g(\bm{x})
		-g(\bm{x},\bm{v})-k\rho_g(\bm{x}).
	\end{align}
\end{subequations} 
By taking the Fourier transform of Eqs.~(\ref{eq_S}) and (\ref{eq_linear}), we obtain
\begin{equation}  
	\hat g(\bm{\lambda},\bm{v})
	=\frac{1-k+\mathrm{i}\frac{F'[0]}{1+d|\bm{\lambda}|^2}\bm{\lambda}\cdot\bm{v}}{1+k\mu+\mathrm{i}k\bm{\lambda}\cdot\bm{v}}\hat\rho_g(\bm{\lambda}).
\end{equation}
Hereafter, $\mathrm{i}$ represents the imaginary unit, $\bm{\lambda}\in\mathbb{R}$ the wave vector, and $\hat a(\bm{\lambda})$ the Fourier transform of the function $a(\bm{x})$ as $\hat a(\bm{\lambda})=\int_{\mathbb{R}} a(\bm{x})e^{-\mathrm{i}\bm{\lambda}\cdot\bm{x}}d\bm{x}$.  
By integrating the above equation as to $\bm{v}$ on $V$, we obtain an equation for $\hat \rho_g(\bm{\lambda})$ as
\begin{subequations}\label{eq_l_rho}
	\begin{gather}
		\hat \rho_g(\bm{\lambda})=\frac{1}{2}\int_{-1}^1
		\frac{1-k+\mathrm{i}\frac{F'[0]\lambda v}{1+d\lambda^2}}{1+k\mu+\mathrm{i}k\lambda v}dv\hat \rho_g(\bm{\lambda}),\\
		\hat \rho_g(\bm{\lambda})=\frac{1}{2}\int_{-1}^1
		\frac{\left(1-k+\mathrm{i}\frac{F'[0]\lambda v}{1+d\lambda^2}\right)
		\left(1+k\mu_1-\mathrm{i}k\lambda(\mu_2+v)\right)}
		{(1+k\mu_1)^2+k^2\lambda^2(\mu_2+v)^2}dv
		\rho_g(\bm{\lambda}),
	\end{gather}
\end{subequations}
where $\mu_1$ is the real part of $\mu$, i.e., $\mu_1=\mathrm{Re}(\mu)$, and $\mu_2$ is defined as $\mu_2=\mathrm{Im}(\mu)/\lambda$.
Here we write $\lambda=|\bm{\lambda}|$.

Thus, in order to obtain a non-trivial solution of $\hat \rho_g(\bm{\lambda})$, the following equations must be simultaneously satisfied, i.e.,
\begin{equation}\label{eq_I1}
	\int_{-1}^1
	\frac{(1-k)(1+k\mu_1)+\frac{F'[0]k\lambda^2v}{1+d\lambda^2}(\mu_2+v)}
	{(1+k\mu_1)^2+k^2\lambda^2(\mu_2+v)^2}dv=2,
\end{equation} 
and
\begin{equation}\label{eq_I2}  
	\int_{-1}^1
	\frac{(1-k)k\lambda(\mu_2+v)-\frac{F'[0]\lambda v}{1+d\lambda^2}(1+k\mu_1)}{(1+k\mu_1)^2+k^2\lambda^2(\mu_2+v)^2}dv=0.
\end{equation} 
Further, Eqs.~(\ref{eq_I1}) and (\ref{eq_I2}) are analytically calculated as, respectively,

\begin{equation}\label{eq_I1x}
	\begin{split}
		\left(\alpha-\frac{\beta}{\xi}\right)
		\left[\arctan\left(\xi(\mu_2+1)\right)
		-\arctan\left(\xi(\mu_2-1)\right)\right]&\\
		+\mu_2 \beta\log\left(\frac{\xi^{-2}+(\mu_2-1)^2}{\xi^{-2}+(\mu_2+1)^2}\right)
		&=2-2\beta,
	\end{split}
\end{equation} 
and
\begin{equation}\label{eq_I2x} 
	\begin{split}
		\mu_2 \beta
		\left[\arctan\left(\xi(\mu_2+1)\right)
		-\arctan\left(\xi(\mu_2-1)\right)\right]&\\
		+\frac{1}{2}\left(\alpha-\frac{\beta}{\xi}\right)
		\log\left(1+\frac{4\mu_2}{\xi^{-2}+(\mu_2-1)^2}\right)
		&=0,
	\end{split}
\end{equation}
where
%
\begin{gather}\label{eq_Isub}
		\alpha=\frac{1-k}{k\lambda},\quad
		\beta=\frac{F'[0]}{k(1+d\lambda^2)},\quad
		\xi=\frac{k\lambda}{1+k\mu_1}.
\end{gather}

Note that Eqs.~(\ref{eq_I1x}) and (\ref{eq_I2x}) are symmetric as to the sign of $\mu_2$ and $\mu_2$=0 always satisfies Eq.~(\ref{eq_I2x}) irrespective of the values of $\alpha$, $\beta$, and $\xi$.
The eigenvalue $\mu\,(=\mu_1+\mathrm{i}\mu_2\lambda)$ is obtained by solving Eqs.~(\ref{eq_I1x}) and (\ref{eq_I2x}) simultaneously, and the sign of the growing rate $\mu_1$ determines the instability of the uniform solution of the kinetic equation, Eq.~(\ref{eq_kinetic}).

Before we consider the instability condition, we first verify that Eq.~(\ref{eq_I1x}) is never satisfied as $\lambda\rightarrow \infty$.
This is explained as follows.

We consider the first term of the left hand side of Eq.~(\ref{eq_I1x}).
In Eq.~(\ref{eq_Isub}), $\alpha$ and $\beta$ vanishes as $\lambda\rightarrow \infty$ while the limiting values of $\xi$ and $\mu_2$ are unknown.
In the case that $\xi$ converges to a finite value or diverges as $\lambda\rightarrow \infty$, the first term of the L.H.S of Eq.~(\ref{eq_I1x}) vanishes because the first factor vanishes while the second factor is bounded.
On the other hand, if $\xi$ converges to zero as $\lambda\rightarrow\infty$, the second factor of the first term of the L.H.S of Eq.~(\ref{eq_I1x}) is estimated as
\begin{equation} 
	\begin{split}
		\left|\arctan(\xi(\mu_2+1))-\arctan(\xi(\mu_2-1))\right|
		=\left|\arctan\left(\frac{2\xi}{1+\xi^2(\mu_2^2-1)}\right)\right|
		\\
		<\left|\arctan\left(\frac{2\xi}{1-\xi^2}\right)\right|
		\rightarrow |2\xi+{\cal O}(\xi^2)|,
	\end{split}
\end{equation}
so that the first term of the L.H.S of Eq.~(\ref{eq_I1x}) vanishes as $\lambda\rightarrow\infty$.
Thus, the first term of the L.H.S of Eq.~(\ref{eq_I1x}) vanishes as $\lambda\rightarrow\infty$ regardless the limiting values of $\xi$ and $\mu_2$.
However, the second term of the L.H.S of Eq.~(\ref{eq_I1x}) is always non-positive while the R.H.S of Eq.~(\ref{eq_I1x}) converges to 2 as $\lambda\rightarrow\infty$.
Thus, Eq.~(\ref{eq_I1x}) cannot be satisfied in the limit of large $\lambda$ whatever the limiting values of $\xi$ and $\mu_2$ take, so that eigenmodes cannot exhibit large oscillations.


\subsection{Instability condition}
We now assume the imaginary part of eigenvalue is zero, i.e., $\mu_2=0$, and consider a linear instability condition, i.e., $\mu_1>0$.
Because $\mu_2=0$ satisfies Eq.~(\ref{eq_I2x}), we only consider Eq.~(\ref{eq_I1x}) with $\mu_2$=0, i.e.,
\begin{equation}\label{eq_I1xx}
	\left(\alpha\xi-\beta\right)
	\frac{\arctan(\xi)}{\xi}
	=1-\beta.
\end{equation}
We can write $\mu_1$ as a function of $\xi$, i.e., 
\begin{equation}\label{e q_mu1x}
	\mu_1=\frac{\lambda}{\xi}-\frac{1}{k},
\end{equation}
so that $\mu_1$ takes a positive value if and only if $0<\xi< k\lambda$.
Thus, the instability condition, i.e., $\mu_1>0$, is equivalent to the condition that Eq.~(\ref{eq_I1xx}) has a solution with $0<\xi< k\lambda$.  

First, we consider Eq.~(\ref{eq_I1xx}) with $\beta=1$, which is obtained at $\lambda=\sqrt{\frac{1}{d}(F'[0]/k-1)}$ when $F'[0]/k\ge1$.
It is immediate that Eq.~(\ref{eq_I1xx}) gives $\mu_1=-1$, so that $\beta=1$ does not hold the instability condition.
(We also note that $\mu_1=-1$ is also obtained as $\lambda$=0.)

For $\beta\ne 1$, we consider the condition that the intersection of the following functions,
\begin{equation}\label{eq_phi}  
	\phi(\xi)=\frac{\arctan(\xi)}{\xi},
\end{equation}
and
\begin{equation}\label{eq_psi}  
	\psi(\xi)=\frac{1-\beta}{\alpha\xi-\beta},
\end{equation}
exists in $0<\xi<k\lambda$.
Note that $\phi(\xi)$ is monotonically decreasing for $\xi\ge 0$ and $\phi(0)=1$ and $\phi(+\infty)=0$, while the behavior of $\psi(\xi)$ depends on the coefficients $\alpha$ and $\beta$.

\begin{figure}[p] 
	\centering
	\includegraphics*[scale=1.0]{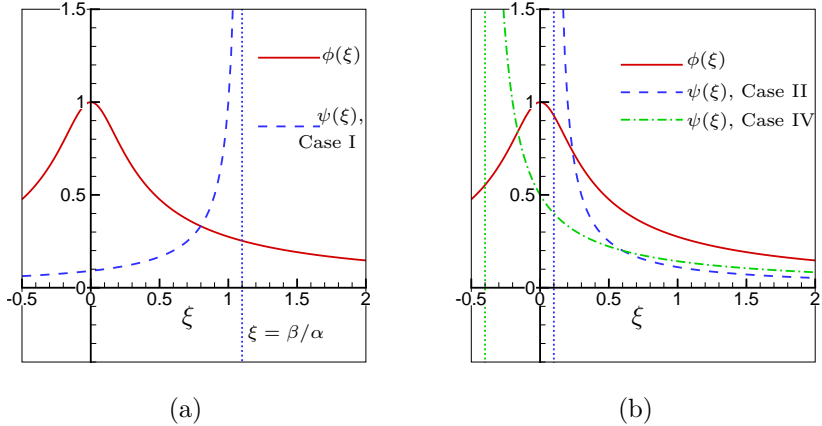}
	\caption{An auxiliary figure for the analysis below Eqs.~(\ref{eq_phi}) and (\ref{eq_psi}). In figure (a), $\alpha=1$ and $\beta=1.1$ is used for Case I. In figure (b), $\alpha=5$ and $\beta=0.5$ is used for Case II and $\alpha=-5$ and $\beta=2$ is used for Case IV. The vertical dotted lines show the asymptotic lines $\xi=\beta/\alpha$ for each case.}\label{fig_aux_analy}
\end{figure}

\noindent\underline{Case I: $\alpha>0$ and $\beta>1$}

The coefficient of $(\xi-\beta/\alpha)^{-1}$ in Eq.~(\ref{eq_psi}) is negative, so that $\psi(\xi)$ takes positive values only when $\xi$ is smaller than the asymptotic line, $\xi=\beta/\alpha$. See also Fig.~\ref{fig_aux_analy}(a).
Thus, we only consider the intersection of $\phi(\xi)$ and $\psi(\xi)$ in $0<\xi<\beta/\alpha$. 
It is found that $\psi(\xi)$ monotonically increases in $0<\xi<k\lambda$ because the asymptotic line, $\xi=\beta/\alpha$, is always larger than $\xi=k\lambda$, because $\frac{\beta}{\alpha}>\frac{1}{\alpha}=\frac{k\lambda}{1-k}>k\lambda$.
It is also found that $\psi(\xi)$ is always smaller than $\phi(\xi)$ at $\xi=0$, because $\psi(0)=1-\frac{1}{\beta}<1$.
Thus, the condition that $\phi(\xi)$ and $\psi(\xi)$ intersect in $0<\xi<k\lambda$ is given by
\begin{equation}\label{eq_fg} 
	\phi(k\lambda)<\psi(k\lambda).
\end{equation}

\noindent\underline{Case II: $\alpha>0$ and $0<\beta<1$}

The coefficient of $(\xi-\beta/\alpha)^{-1}$ in Eq.~(\ref{eq_psi}) is positive, so that $\psi(\xi)$ takes positive values only when $\xi>\beta/\alpha$.
It is obvious that if the asymptotic line $\xi=\beta/\alpha$ is larger than or equal to $\xi=k\lambda$, $\phi(\xi)$ and $\psi(\xi)$ do not intersect in $0<\xi<k\lambda$.
See also Fig.~\ref{fig_aux_analy}(b).
Thus, we only consider the case $\beta/\alpha<k\lambda$, in which $\psi(\xi)$ is positive and monotonically decreases in $\beta/\alpha<\xi<k\lambda$.
However, $\psi(\xi)$ is always larger than the unity in this range.
This is explained as follows.

From the condition $\beta/\alpha<k\lambda$ and Eq.~(\ref{eq_Isub}), we obtain
\begin{equation} 
	0<\frac{k}{1-\beta}<1.
\end{equation}
Since $\psi(\xi)$ monotonically decreases in $\beta/\alpha<\xi<k\lambda$, we get 
\begin{equation} 
	\psi(\xi)>\psi(k\lambda)=\frac{1}{1-\frac{k}{1-\beta}}>1.
\end{equation}
Thus, we can conclude that $\phi(\xi)$ and $\psi(\xi)$ do not intersect in $0<\xi<k\lambda$ in this case. 

\noindent\underline{Case III: $\alpha=0$, i.e., $k=1$.}

In this case, $\psi(\xi)=1-\frac{1}{\beta}$ is constant, so that the condition that $\phi(\xi)$ and $\psi(\xi)$ intersect in $0<\xi<k\lambda$ is given by Eq.~(\ref{eq_fg}). 

\noindent\underline{Case IV: $\alpha<0$ and $\beta>1$.}

The coefficient of $(\xi-\beta/\alpha)^{-1}$ in Eq.~(\ref{eq_psi}) is positive and the asymptotic line $\xi=\beta/\alpha$ is in the negative region,  so that $\psi(\xi)$ takes finite positive values and monotonically decreases in $0<\xi<k\lambda$.
See also Fig.~\ref{fig_aux_analy}(b).
In addition, $\psi(\xi)$ is smaller than $\phi(\xi)$ at $\xi=0$, because $\psi(0)=1-1/\beta<1$.
Thus, the condition that $\phi(\xi)$ and $\psi(\xi)$ intersect in $0<\xi<k\lambda$ is given by Eq.~(\ref{eq_fg}).

\noindent\underline{Case V: $\alpha<0$ and $0<\beta<1$.}

The coefficient of $(\xi-\beta/\alpha)^{-1}$ in Eq.~(\ref{eq_psi}) is negative and the asymptotic line $\xi=\beta/\alpha$ is in the negative region, so that $\psi(\xi)$ is always negative in $0<\xi<k\lambda$.
Thus, $\phi(\xi)$ and $\psi(\xi)$ do not intersect in $0<\xi<k\lambda$.

In summary, the condition that $\phi(\xi)$ and $\psi(\xi)$ intersect in $0<\xi<k\lambda$ is solely given by Eq.~(\ref{eq_fg}). Thus, the instability condition of the perturbation with mode $\lambda$ for Eqs.~(\ref{eq_kinetic})--(\ref{eq_P}) is written as
\begin{equation}\label{ineq_kinetic_instability}
	\frac{F'[0]}{k}>\left[
		1+\frac{k}{\frac{k\lambda}{\arctan(k\lambda)}-1}
	\right]
	(1+d\lambda^2).
\end{equation}
We also remark that the auxiliary condition $\beta>1$, i.e., $0<\lambda<\sqrt{(F'[0]/k-1)/d}$, is automatically satisfied under the above condition.
Thus, when $F'[0]/k$ exceeds the minimum of the right-hand side of Eq.~(\ref{ineq_kinetic_instability}) with variation in  $0<\lambda<\sqrt{(F'[0]/k-1)/d}$, the stationary homogeneous state in the population density becomes linearly unstable.
Furthermore, stationary periodic patterns are generated because no unstable eigenmodes exist in the limit of large $\lambda$.

\subsection{Continuum limit} 
We introduce a small parameter $\varepsilon$ and scale $k$, $x$, and $d$ in Eqs.~(\ref{eq_kinetic}) and (\ref{eq_S}) as follows,
\begin{gather}\label{eq_scale}
	k=\varepsilon^2,\quad x=\varepsilon \hat x,
	\quad d=\varepsilon^2\hat d.
\end{gather}
We also suppose the gradient of the response function is scaled as
\begin{gather}
	F'[0]=\varepsilon^2\hat F'[0].
\end{gather}
Then, by the asymptotic analysis of small $\varepsilon$ for Eqs.~(\ref{eq_kinetic}) and (\ref{eq_S}), we obtain the following flux-limited Keller-Segel equation in the limit $\varepsilon\rightarrow 0$,
\begin{equation}\label{eq_rho0}
	\partial_t\rho+\partial_{\hat x}\left(U[\log S]\rho\right)
	=\frac{1}{3}\partial_{\hat x\hat x}\rho+P[\rho]\rho,
\end{equation} 
\begin{equation}\label{eq_S0}
	-\hat d\partial_{\hat x\hat x}S+S=\rho,
\end{equation}
where $U[\log  S]$ is defined as
\begin{equation}\label{eq_U}
	U[\log S]=\int_0^1v \hat F[v\partial_{\hat x}\log S]dv.
\end{equation}
We remark that the flux $U[\log S]$ is bounded as $|U[\log S]|\le \frac{\chi}{2}$ because of the bounded signal response Eq.~(\ref{cond_F}c).

When we carry out the linear instability analysis for Eqs.~(\ref{eq_rho0})--(\ref{eq_U}) as in the previous subsection, we obtain the growth rate
\begin{equation}\label{eq_mu1_continuum}
	\mu_1=-1+\frac{\hat F'[0]\hat \lambda^2}{3(1+\hat d\hat \lambda^2)}-\frac{1}{3}\hat \lambda^2,
\end{equation}
where $\hat \lambda$ is the Fourier variable as to $\hat x$, i.e., $\hat \lambda=\varepsilon \lambda$.
Eq.~(\ref{eq_mu1_continuum}) achieves its  maximum at $\hat \lambda^2=(\sqrt{\hat F'[0]}-1)/\hat d$, and the sign of the maximum value of the growth rate $\mu_1$ determines the linear stability.
Thus, the instability condition for Eqs.~(\ref{eq_rho0})--(\ref{eq_U}) is written as
\begin{equation}\label{continuum_instability}
	\hat F'[0]>(1+\sqrt{3\hat d})^2.
\end{equation}
Remarkably the above equation is consistent with the linear stability condition of the Keller-Segel system obtained earlier by Nadin, {\it et. al.,} Ref. {\cite{art:07NPR}.
We also note that Eq.~(\ref{continuum_instability}) can be obtained by using Eq.~(\ref{kinetic_instability}) with the same scaling as Eq.~(\ref{eq_scale}).
This means Eq.~(\ref{kinetic_instability}) is not only sufficient but also necessary condition for linear instability in the continuum limit.

\section{Numerical analysis}
The numerical simulations are performed for Eqs.~(\ref{eq_kinetic})--(\ref{eq_rho}) with the uniform initial density, i.e., $f(0,x,v)=S(0,x)=1$, and periodic boundary condition in the one-dimensional interval $x=[0,L]$, i.e., $f(t,0,v)=f(t,L,v)$ and $S(t,0)=S(t,L)$.
The kinetic equation Eq.~(\ref{eq_kinetic}) is solved by the Monte Carlo (MC) method, which has been recently developed in Ref. \cite{art:17Y}, coupled with the finite volume (FV) scheme for Eq.~(\ref{eq_S}).
The details of the MC method and FV scheme is presented in Sec. \ref{sec_mcmethod}.
In the MC simulations, the one-dimensional interval $L$ is set as $L$=100 and divided into 2000 cubic lattice boxes with a side length $\Delta x$=0.05.
The 1$\times 10^7$ simulation particles are used initially as a total in the whole lattice boxes, and they are distributed randomly into each lattice box. See also Fig.~\ref{fig_geom}.
The time-step size $\Delta t$ is set as $\Delta t$=5$\times 10^{-3}$.
In the previous study, it has been verified that these simulation parameters produce accurate numerical solutions.
Eq.~(\ref{eq_S}) is discretized using the FV scheme on the uniform lattice mesh system with a mesh interval $\Delta x$ and solved implicitly at each time step.

For the population-growth rate $P[\rho]$ and response function $F[X]$, which satisfies Eq.~(\ref{eq_P}) and Eq.~(\ref{cond_F}), respectively, we consider
\begin{equation}\label{eq_logistic}
	P[\rho]=1-\rho,
\end{equation}
and
\begin{equation}\label{eq_tanh}
	F[X]=\chi\tanh\left(\frac{X}{\delta}\right).
\end{equation}
Here, $\chi$ and $\delta^{-1}$ represent the amplitude of modulation and stiffness of the response function, respectively.
Note that for Eq.~(\ref{eq_tanh}), $F'[0]$ is given as $F'[0]=\chi/\delta$.
The values of parameters $\chi$, $\delta$, and $d$ are set so as to correspond to the marks A, B, C in Fig.~\ref{fig_dki} with four values of  $k$, i.e., $k$=0.1, 1, 2, and 10.
More specifically, we set the parameters as $(\chi/\sqrt{k}, \sqrt{k}\delta, d/k)$=(0.5, 0.05, 1) for A, (0.5, 0.0625, 1) for B, (0.5, 0.0625, 0.7) for C, and (0.5, 0.1, 1) for D. 
We perform the MC simulations for the parameters listed in Table \ref{t_para} and investigate the compatibility and sharpness of the kinetic instability condition Eq.~(\ref{kinetic_instability}) numerically.
Table \ref{t_para} also shows the prediction of the linear instability by Eq.~(\ref{kinetic_instability}).
For example, the parameter set C with $k=1.0$ is slightly above the critical line in Fig.~\ref{fig_dki}, so that the kinetic instability condition Eq.~(\ref{kinetic_instability}) affirms the occurrence of periodic pattern formation.
However, for the parameter set C with $k=2.0$, which is slightly below (but very close to) the critical line in Fig.~\ref{fig_dki}, we cannot confirm neither the linear instability nor homogeneous state theoretically because Eq.~(\ref{kinetic_instability}) is a sufficient condition of the linear instability.
However, the MC simulation can numerically demonstrate how sharply the kinetic instability condition Eq.~(\ref{kinetic_instability}) can predict the pattern formations.
\begin{table}[p] 
	\centering
	\begin{tabular}{ccccc}
		\hline\hline
		$k$& A(1, 0.5, 0.05)  & B(1, 0.5, 0.0625) & C(0.7, 0.5, 0.0625) & D(1, 0.5, 0.1)\\
		\hline
		0.1& -- & $\blacksquare$ & -- &--\\
		1.0& $\blacksquare$  & $\square$ & $\blacksquare$ & $\square$\\
		2.0& $\square$  & -- & $\square$ &--\\
		\hline\hline
	\end{tabular}
	\caption{The parameter sets with which the Monte Carlo simulations are performed.
		The parameter sets A, B, C, and D corresponds to the marks A, B, C, and D in Fig.~\ref{fig_dki}, respectively.
		The values of the round brackets show the values of $d/k$, $\chi/\sqrt{k}$, and $\sqrt{k}\delta$, i.e., ($d/k$, $\chi/\sqrt{k}$, $\sqrt{k}\delta$).
		The black squares satisfy the kinetic instability condition Eq.~(\ref{kinetic_instability}) while white squares do not satisfy it.}\label{t_para}

	\end{table}
	\begin{figure}[p]
		\centering
		\includegraphics*[scale=0.9]{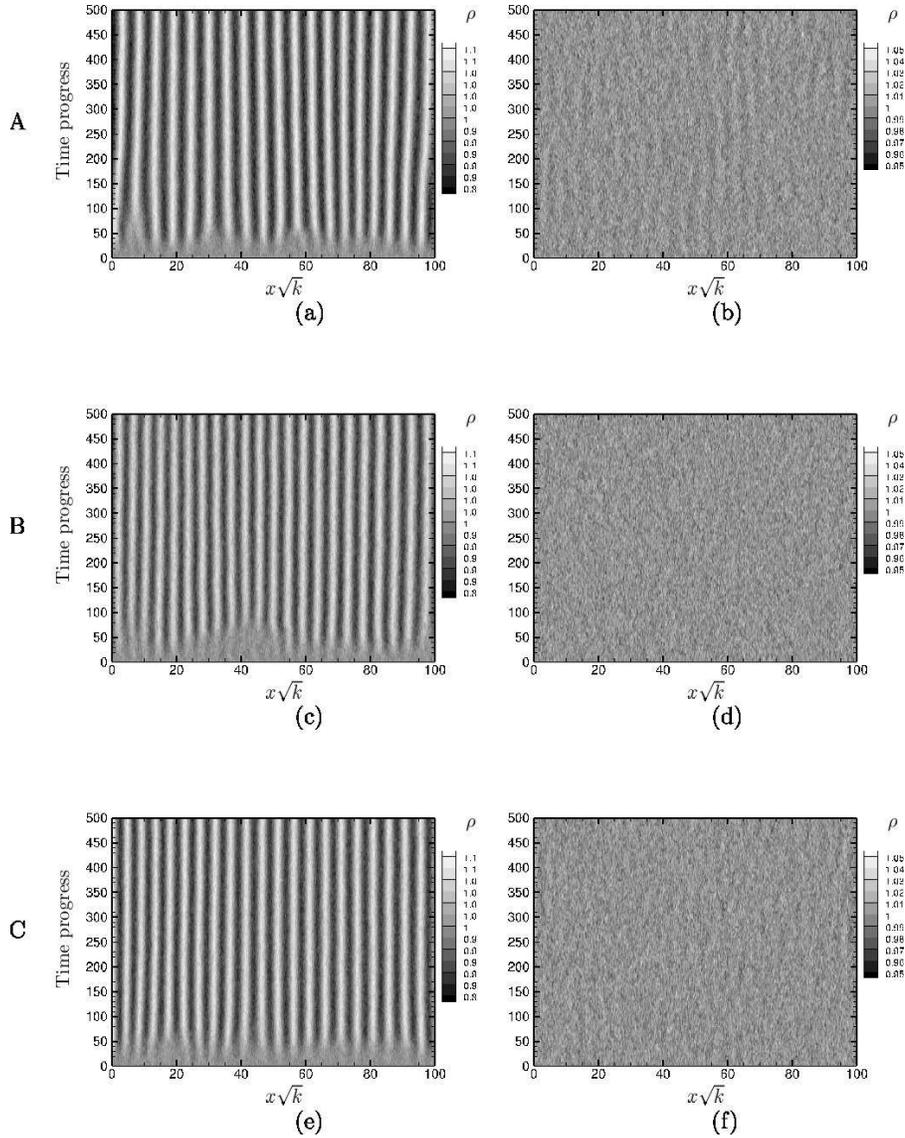}
		\caption{
			The time progress of the population densities.
			Figure (a) and (b) correspond to the parameter set A with $k$=1.0 and 2.0, respectively, Figure (c) and (d) to the parameter set B with $k$=0.1 and 1.0, respectively, and Figure (e) and (f) to the parameter set C with $k$=1.0 and 2.0, respectively. See also Table \ref{t_para}.
		}\label{fig_2dmap}
	\end{figure}
	\begin{figure}[p]
		\centering
		\includegraphics*[scale=0.8]{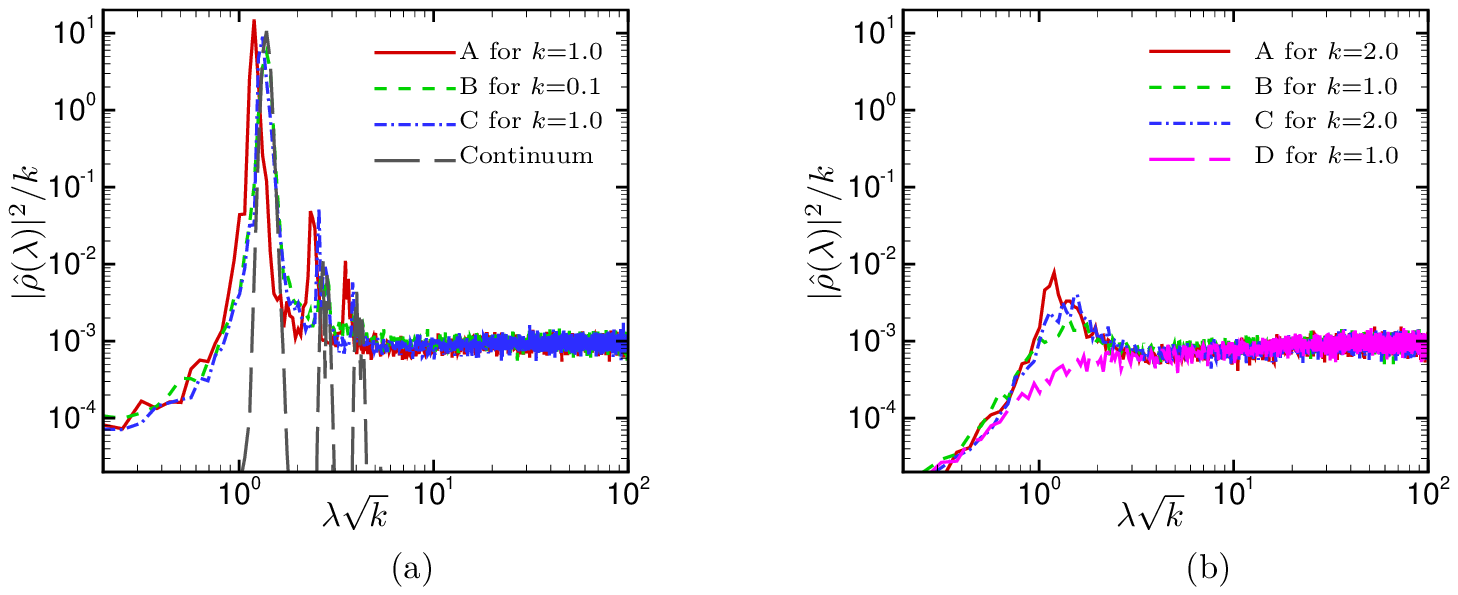}
		\caption{The power spectra of the Fourier transform of population density profiles, $|\hat \rho(\lambda)|^2/k$. Figure (a) shows the results for the parameter sets A with $k=1.0$, B with $k=0.1$, and C with $k=1.0$ (the black squares $\blacksquare$) in Table \ref{t_para} and figure (b) shows those for A with $k=2.0$, B with $k=1.0$, C with $k=2.0$, and D with $k=1$ (the white squares $\square$) in Table \ref{t_para}.
		In the figure (a), the result of the continuum equations Eqs.~(\ref{eq_rho0})--(\ref{eq_U}) with the parameter set B in Table \ref{t_para} (the dotted line) is also included.
	}\label{fig_turing_pow}
\end{figure}

Figure \ref{fig_2dmap} shows the time progress of the population density $\rho$.
It is obviously seen that the results for the black squares in Table \ref{t_para} exhibit the stationary periodic patterns after some transient period.
On the other hand, those for the white squares in Table \ref{t_para} do not show any distinct patterns.
In order to quantify the patterns, we also calculate the power spectra of the Fourier transforms of population density profiles between $t$=[400,500].
Figure \ref{fig_turing_pow} shows the results of the power spectra.
Here, we take the averages of the snap shots of power spectra obtained at $t$=[400,500] with a time interval 4.
The power spectrum for the continuum equations Eqs.~(\ref{eq_rho0})--(\ref{eq_U}) is calculated from the snap shot of the population density at $t=200$.
For the parameter sets shown as the black squares $\blacksquare$ in Table \ref{t_para} (See Fig.~\ref{fig_turing_pow}(a)), steep peaks are observed and the first peaks appear in $0<\hat \lambda<\sqrt{(\hat F'[0]-1)/\hat d}$ for each parameter set.
The second and third peaks also appear as the non-linear effects although they are much smaller than the first peaks.
The power spectra decreases as the wave number $\hat \lambda$ decreases from the first-peak position, while they neither grow nor damp at the large wave number, so that a plateau regime appears at the large wave numbers.
The peaks of the power spectrum for the continuum equations with the parameter set B coincide with those for the parameter set B with $k=0.1$.
However, no plateau regime appears for the continuum equations.
This result confirms that there are no eigenmodes of the linearized kinetic equation.

For the parameter sets shown as the white squares $\square$ in Table \ref{t_para} (See Fig.~\ref{fig_turing_pow}(b)), the behaviors of power spectra are similar to those for the black squares in Table \ref{t_para}, i.e., Fig.~\ref{fig_turing_pow}(a), except the peak behaviors.
In Fig.~\ref{fig_turing_pow}(b), we cannot see the steep first peaks as seen in Fig.~\ref{fig_turing_pow}(a).
The small peaks of A, B, and C in Fig.~\ref{fig_turing_pow}(b) are even smaller than the second peaks appearing in the figure (a).
For the parameter set C, we cannot see any peaks.
Thus, we cannot observe the stationary periodic patterns evidently in Fig.~\ref{fig_2dmap} for the parameter sets shown as white squares in Table \ref{t_para}.

The parameter set A and B are very close to the critical lines for $k$=2.0 and $k$=0.1 in Fig.~\ref{fig_dki}, respectively.
(The parameter set A is only slightly lower and B is only slightly upper than each critical line.)
However, the result of the parameter set B shows the stationary periodic pattern evidently but the result of the parameter set A does not show it.
These numerical results demonstrate that the critical lines in the instability diagram Fig.~\ref{fig_dki} can predict sharply the occurrence of the linear instability and pattern formations.

\begin{figure}[h]
	\centering
	\includegraphics*[scale=0.8]{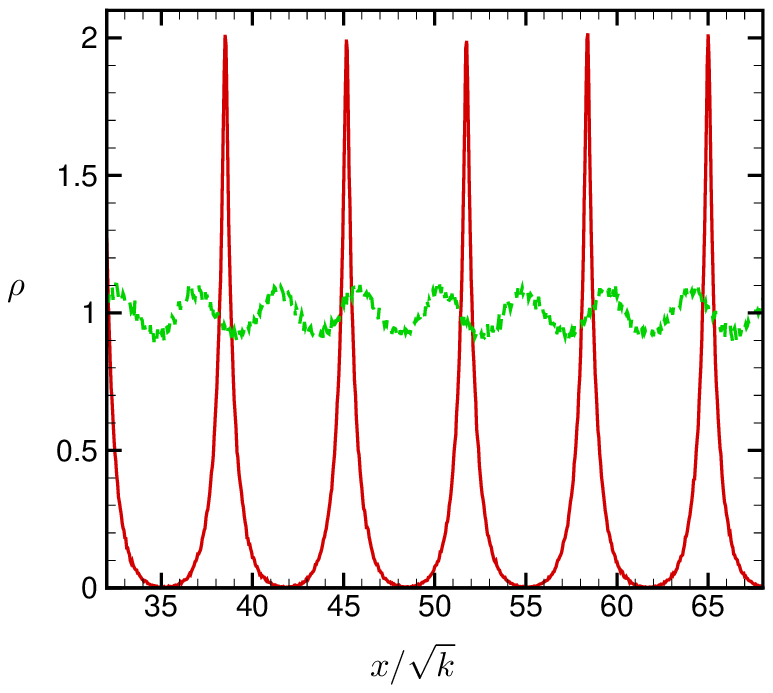}
	\caption{
		The effect of modulation amplitude $\chi$ on the instability pattern profile.
		The dotted line shows the result obtained by the parameter set B with $k=0.1$ in Table \ref{t_para} (which is also shown in Fig. \ref{fig_2dmap} (c)).
		The solid line shows the result obtained when the modulation amplitude is about four-times larger than the parameter set B with $k=0.1$, i.e., $\chi/\sqrt{k}=2.06$.
	}\label{fig_spik}
\end{figure}
Finally, we show the effect of modulation amplitude $\chi$ on the instability pattern profile.
Figure \ref{fig_spik} shows the pattern profile for the parameter set B with $k=0.1$ in Table \ref{t_para} and that obtained when the modulation amplitude $\chi$ is about four times larger than the parameter set B for $k=0.1$.
The pattern profile obtained for the former parameter is periodic oscillation around the initial uniform state $\rho=1$.
The formation of periodic oscillatory patter applies to all of the instability patterns obtained in the parameter sets in Table \ref{t_para}, where $\chi$ is fixed as $\chi/\sqrt{k}=0.5$.
However, in the case of a large modulation amplitude, i.e., $\chi/\sqrt{k}=2.06$, the population of bacteria becomes localized due to a strong chemotactic response so that it forms periodical bounded spikes.
The boundedness in instability pattern formation stems from the flux-limited property in the non-linear stiff response function.
Incidentally, the boundedness property is also inherent in the flux-limited Keller-Segel equation which is obtained by the asymptotic analysis of the kinetic chemotaxis equation.
The variety of solution types is observed in the flux-limited Keller-Segel equation, which will be addressed in our forthcoming paper.

\section{Monte Carlo method}\label{sec_mcmethod}
\begin{figure}[p]
	\centering
	\includegraphics*[scale=1.0]{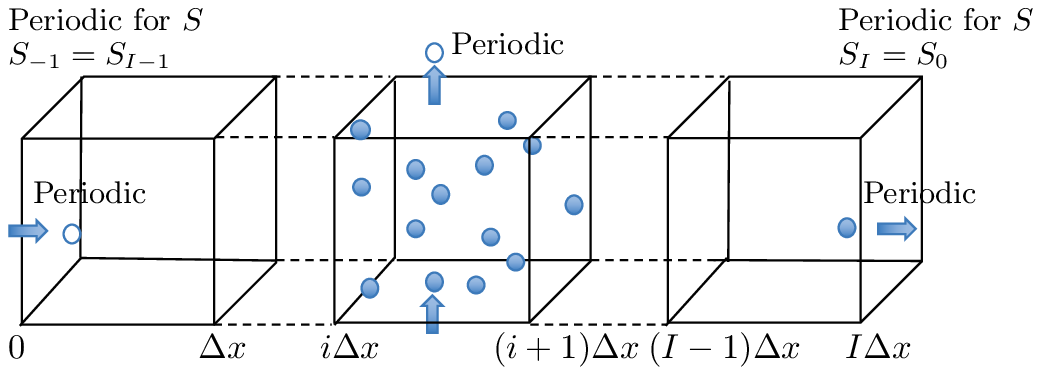}
	\caption{The geometry in the Monte Carlo method. The concentrations of chemoattractant are calculated by a finite-volume method on the uniform cubic lattice boxes. The Monte Carlo particles are distributed in each lattice boxes. The periodic boundary conditions are considered both for the MC particles and chemoattractant in this paper.}\label{fig_geom}
\end{figure}
The motions of the chemotactic bacteria are simulated by using the Monte Carlo particles which follows the process described by the kinetic chemotaxis equation Eq.~(\ref{eq_kinetic}) coupled with the reaction-diffusion equation for the chemoattractant Eq.~(\ref{eq_S}).
The one-dimensional space interval $x\in[0,L]$ is divided into $I$ cubic lattices with a uniform side length $\Delta x$, i.e., $L=I\Delta x$.
See also Fig.~\ref{fig_geom}.
The reaction-diffusion equation Eq.~(\ref{eq_S}) is implicitly solved on the uniform lattice system by using a finite volume scheme as
\begin{equation}\label{fv_S}
	-\frac{d}{\Delta x^2}(S_{i+1}^n-2S_i^n+S_{i-1}^n)+S_i^n=\rho^n_i,
	\quad(i=0,\cdots,I-1),
\end{equation}
where $S_i^n$ and $\rho_i^n$ are the concentration of chemical attractant and population density of bacteria in the $i$th lattice site $[i\Delta x,(i+1)\Delta x]\times \Delta x^2$ at a time $t=n\Delta t$, respectively.
Here we also set $S^n_{-1}=S^n_{I-1}$ and $S^n_{I}=S^n_0$ according to the periodic condition.
Hereafter the superscript represents the time-step number, the subscript without parenthesis represents the lattice-site number, and the subscript in parenthesis represents the index of each MC particle. 
The population density of bacteria $\rho^n_i$ is calculated from the number of the MC particles involved in the $i$th lattice site, $M_i^n$, as
\begin{equation}\label{eq_rho_i_n}  
	\rho^n_i=M^n_i/M,
\end{equation}
where $M$ is the number of MC particles involved in one lattice site in the reference state, i.e, $M=N/I$, where $N$ is the total number of MC particles in the initial state.

The MC simulation is conducted using the following steps. Hereafter, the position and velocity of the $l$th particle are expressed as ${\bm r}_{(l)}^n$ and ${\bm v}_{(l)}^n$, respectively.

\begin{enumerate}
		\setcounter{enumi}{-1}
	\item At $\hat t=0$, MC particles are distributed according to the initial density.
		In each lattice site, MC particles are distributed uniformly at random positions and their velocities ${\bm v}$ are determined by the probability density $f_i^0({\bm v})/(4\pi\rho_i^0)$.
	\item 
		Particles move with their velocities for a duration $\Delta \hat t$:
		\begin{equation}
			{\bm r}_{(l)}^{n+1}={\bm r}_{(l)}^n+{\bm v}_{(l)}^n\Delta \hat t\quad (l=1,\cdots,N^n),
		\end{equation}
		where $N^n$ is the total number of simulation particles at a time step $n$, i.e., $N^n=\sum_{i=0}^{I-1}M_i^n$.
		The particles that move beyond the boundaries are inserted at the opposite boundaries according to the periodic boundary conditions.
	\item
		At each lattice site, the macroscopic population densities $\rho^{n+1}_i$ and concentrations of chemical cues $S^{n+1}_i$ ($i=0,\cdots,I-1$) are calculated by Eqs.~(\ref{eq_rho_i_n}) and (\ref{fv_S}), respectively. 
	\item  
		The tumbling of each particle is calculated using the scattering kernel in Eq.~(\ref{eq_kinetic}).
		The tumbling of the $l$th particle may occur with a probability
		\begin{equation}
			\frac{\Delta t}{k}
			K[D_t\log S^{n+1}_{(l)}],
		\end{equation}
		where $D_t \log S_{(l)}$ represents the temporal variation of the chemical cue experienced by the $l$th MC particle along the pathway, and is defined by the following forward difference,
		\begin{equation}\label{eq_DtlogS_l}
			D_t\log S^n_{(l)}=\frac{\log S(t^n,{\bm r}^n_{(l)})
			-\log S(t^{n}-\Delta t,{\bm r}^{n}_{(l)}-{\bm v}^n_{(l)}\Delta t)}{\Delta t}.
		\end{equation}
		The local concentration of chemical cue at the position of the $l$th MC particle is calculated by using linear interpolation between the neighboring lattice sites, i.e., 
		\begin{equation}\label{eq_interp}
			\log S(t^n,{\bm r}^n_{(l)})=\log S^n_i + 
			\frac{\log S^n_{i+1}-\log S^n_{i-1}}{2\Delta x}
			({r_x}^n_{(l)}-{x}_{i+\frac{1}{2}}).
		\end{equation}
		This generates a chemoattractant gradient and MC particles that stay at the same lattice site after a single time step can sense the chemoattractant gradient along their pathways.

		For the particle that is judged to tumble, say the $l_t$th particle, a new velocity after the tumbling, ${\bm v}^{n+1}_{(l_t)}$, is determined randomly as,
		\begin{equation}\label{random_vel}
			v_x=1-2U_1,\quad v_y=\sqrt{1-v_x^2}\cos(2\pi U_2), \quad v_z=\sqrt{1-v_x^2}\sin(2\pi U_2).
		\end{equation}
		Here $U_1$ and $U_2$ are the uniform random variables between 0 and 1.
	\item
		The divisions/deaths are judged for all MC particles.
		The division (or death) occurs with a probability $|P[\rho^n_i]|\Delta t$, if $|P[\rho_i^n]|$ is positive (or negative), where $\rho_i^n$ is the local population density at the lattice site where each MC particle is involved.
		For a particle that is judged to undergo division, e.g., the $l$th particle, a new particle with the same velocity ${\bm v}_{(l)}$ is created at a random position within the same lattice site.
		The numbers of MC particles involved in each lattice site are counted, $M_i^{n+1}$ ($i=0,\cdots,I-1$), and the total number of simulation particles is updated as $N^{n+1}$.
	\item
		Return to step 1 with the obtained ${{\bm r}}_{(l)}$, ${{\bm e}}_{(l)}$, $S_{(l)}$ ($l$=1,$\cdots$,$M$) at the new time step.
\end{enumerate}

\noindent\underline{Weak formulation}
The overall procedure corresponds to the first-order time difference equation of the kinetic chemotaxis equation in the weak formulation.
We consider the following moment equation,
\begin{align}\label{eq_kinetic_dt}
	&<\Phi({\bm x},{\bm v}),f(t+\Delta t,{\bm x},{\bm v})>=
	<\Phi({\bm x},\bm{v}),f(t,\bm{x},\bm{v})-\Delta t \bm{v}\cdot\nabla f(t,\bm{x},\bm{v})> \nonumber \\
	&+<\Phi(\bm{x},\bm{v}),\frac{1}{4\pi k}\int_V K[D_t \log S|_{\bm{v}'}]f(t,\bm{x},\bm{v}')d\Omega(\bm{v}')-\frac{1}{k}K[D_t \log S|_{\bm{v}}]f(t,\bm{x},\bm{v})>\Delta t \nonumber \\
	&+<\Phi(\bm{x},\bm{v}),P[\rho(t,\bm{x})]f(t,\bm{x},\bm{v})>\Delta t,
\end{align}
where  $\Phi(\bm{x},\bm{v})$ is an arbitrary smooth function which vanishes outside the computational domain on~$\bm{x}$.

Here $<\quad,\quad>$ defines the integration of the arbitrary functions $a(\bm{x},\bm{v})$ and $b(\bm{x},\bm{v})$ as
\begin{equation}
	<a(\bm{x},\bm{v}),b(\bm{x},\bm{v})>=\frac{1}{4\pi(\Delta x)^3}\int_{V,\mathbb{R}}
	a(\bm{x},\bm{v})b(\bm{x},\bm{v})d\bm{x}d\Omega(\bm{v}).
\end{equation}
We consider the functions $f^A$, $f^B$, and $f^C$ which are determined, respectively, as
\begin{subequations}
	\begin{equation}\label{eq_fA}
		<\Phi,f^A>=
		<\Phi({\bm x},\bm{v}),f(t,\bm{x},\bm{v})-\Delta t \bm{v}\cdot\nabla f(t,\bm{x},\bm{v})>,
	\end{equation}
	\begin{equation}\label{eq_fB}
		\begin{split}
			<\Phi,f^B>=
			<\Phi(\bm{x},\bm{v}),\frac{1}{4\pi k}\int_V K[D_t \log S|_{\bm{v}'}]f(t,\bm{x},\bm{v}')d\Omega(\bm{v}')\\
			-\frac{1}{k}K[D_t \log S|_{\bm{v}}]f(t,\bm{x},\bm{v})>\Delta t,
		\end{split}
	\end{equation}
	\begin{equation}\label{eq_fC}
		<\Phi,f^C>=<\Phi(\bm{x},\bm{v}),P[\rho(t,\bm{x})]f(t,\bm{x},\bm{v})>\Delta t.
	\end{equation}
\end{subequations}
Then, $f(t+\Delta t,\bm{x},\bm{v})$ is obtained by the sum of three functions:
\begin{equation}\label{sum_fABC}
	f(t+\Delta t,\bm{x},\bm{v})=f^A+f^B+f^C.
\end{equation}

In the MC method, the microscopic population density $f$ at the $n$th time step is approximated as
\begin{equation}\label{eq_distf} 
	f^n(\bm{x},\bm{v})=\frac{4\pi(\Delta x)^3}{M}\sum_{l=1}^{N^n}
	\delta(\bm{x}-\bm{r}_{(l)}^n)\delta(\bm{v}-\bm{v}_{(l)}^n).
\end{equation}
By substituting Eq.~(\ref{eq_distf}) into Eq.~(\ref{eq_fA}) we obtain the following equation:
\begin{equation}
	\begin{split}
		<\Phi,f^A>&=\frac{1}{M}\sum_{l=1}^{N^n} \Phi(\bm{r}_{(l)}^n,\bm{v}_{(l)}^n)+\Delta t \bm{v}\cdot
		\nabla \Phi(\bm{r}_{(l)}^n,\bm{v}_{(l)}^n)\\
		&=\frac{1}{M}\sum_{l=1}^{N^n} \Phi(\bm{r}_{(l)}^n+\bm{v}_{(l)}^n\Delta t,\bm{v}_{(l)}^n)+{\cal O}(\Delta t^2).
	\end{split}
\end{equation}
It is easily seen that the function $f^A$ corresponds to the distribution that is obtained by the moving process, say ${\cal S}_1$ in the MC method, i.e., 
\begin{equation}\label{transS1}
	{\cal S}_1\{f^n\}=f^A+{\cal O}(\Delta t^2), 
\end{equation}
where ${\cal S}_1\{\quad\}$ represents the operator of the process 1 in the MC method. 

The tumbling process and division/death process in the MC method, say ${\cal S}_3$ and ${\cal S}_4$ respectively, are performed independently in each lattice site.
Thus, we consider the processes in a fixed lattice site, say the $i$th lattice site, with using a test function written as $\Phi(\bm{x},\bm{v})=\phi_i(\bm{x})\psi(\bm{v})$, where $\phi_i(\bm{x})$ is an arbitrary smooth function which vanishes outside the $i$th lattice site.
Then, by substituting Eq.~(\ref{eq_distf}) into Eqs.~(\ref{eq_fB}) and (\ref{eq_fC}) we obtain the following equations:
\begin{align}\label{eq_momentfB}
	\begin{split}
		<\phi_i(\bm{x})\psi(\bm{v}),f^B>=
		\frac{\Delta t}{kM}\sum_{m=1}^{M^n_i}
		K\left[D_t\log S_{(m)}^n\right]
		\phi_i(\bm{r}^n_{(m)})(\bar \psi-\psi(\bm{v}^n_{(m)}))\\
		+{\cal O}(\Delta t^2),
	\end{split}
\end{align}
\begin{align}\label{eq_momentfC}
	<\phi_i(\bm{x})\psi(\bm{v}),f^C>&=
	\frac{\Delta t}{M}\sum_{m=1}^{M^n_i}P[\rho(\bm{r}_{(m)}^n)]\phi_i(\bm{r}_{(m)}^n)\psi(\bm{v}_{(m)}^n),\nonumber\\
	&=\frac{\Delta t P[\rho_i^n]}{M}\sum_{m=1}^{M^n_i}\bar \phi_i \psi(\bm{v}_{(m)}^n)+{\cal O}(\Delta t \Delta x^2),
\end{align}
where the subscript $m$ ($m=1,\cdots,M_i^n$) counts the MC particles involved in the $i$th lattice site at the $n$th time step. 
Here $D_t\log S_{(m)}^n$ is defined in Eq.~(\ref{eq_DtlogS_l}),
$\bar \psi=\frac{1}{4\pi}\int \psi(\bm{v})d\Omega(\bm{v})$,
and $\bar \phi_i=\phi(\bm{x}_{i+\frac{1}{2}})$, where $\bm{x}_{i+\frac{1}{2}}$ represents the center of the $i$th lattice site.  
In the derivation of Eq.~(\ref{eq_momentfB}), we approximate the tumbling kernel as
\begin{equation}
	K[D_t\log S|_{\bm{v}}]=K\left[
		\frac{\log S(t,\bm{x})-\log S(t-\Delta t,\bm{x}-\bm{v}\Delta t)}
	{\Delta t}\right]+{\cal{O}}(\Delta t).
\end{equation}
We note that the second equality in Eq. (\ref{eq_momentfC}) is obtained under the assumption of the uniform distribution of large number of particles in each lattice site.

In the tumbling process ${\cal S}_3$, the particle which creates the tumbling, say the $m_t$th particle, changes its velocity randomly as $\bm{v}^n_{m_t}\rightarrow \bm{v}^{n+1}_{m_t}=\bm{u}_t$, where the random velocity $\bm{u}_t$ is given by Eq.~(\ref{random_vel}).
We now suppose that $M_t$ particles make the tumbling in the $i$th lattice site, then the microscopic population density $f^n$ in the $i$th lattice site changes as
\begin{equation}\label{eq_distS3}
	{\cal S}_3\{f^n\}=\frac{4\pi(\Delta x)^3}{M}
	\left\{
		\sum_{\substack{m=1\\ m \ne m_t}}^{M_i^n}
		\delta(\bm{x}-\bm{r}_{(m)}^n)\delta(\bm{v}-\bm{v}_{(m)}^n)
		+
		\sum_{t=1}^{M_t}
		\delta(\bm{x}-\bm{r}_{(m_t)}^n)\delta(\bm{v}-\bm{u}_t)
	\right\},
\end{equation}
and the moment is written as
\begin{equation}\label{eq_momentS3}
	\begin{split}
		&<\phi_i(\bm{x})\psi(\bm{v}),{\cal S}_3\{f^n\}>\\
		&=\frac{1}{M}\left\{
			\sum_{m=1}^{M_i^n}\phi_i(\bm{r}_{(m)}^n)\psi(\bm{v}_{(m)}^n)
			+\sum_{t=1}^{M_t}\phi_i(\bm{r}_{(m_t)}^n)
			\left[\psi(\bm{u}_{t})-\psi(\bm{v}_{(m_t)}^n)\right]
		\right\}.
	\end{split}
\end{equation}
We introduce the stochastic variables $Z_m$ and $Z$ which are defined as
\begin{equation}\label{eq_Zm}
	Z_m=\left\{
		\begin{array}{cl}
			\phi_i(\bm{r}_{(m)})\left[\psi(\bm{u})-\psi(\bm{v}_{(m)})\right],&(\mbox{if tumbling}),\\
			0,&(\mbox{otherwise}),
		\end{array}
		\right.
	\end{equation}
	and
	\begin{equation}\label{eq_Z}
		Z=\frac{1}{M}\sum_{m=1}^{M_i^n}Z_m.
	\end{equation}
	If the tumbling occurs for $M_t$ particles in the $i$th lattice site, the realized value of $Z$ is written as
	\begin{equation}\label{eq_realZ}
		\hat Z=\frac{1}{M}\sum_{t=1}^{M_t}\phi(\bm{r}_{m_t})
		\left[
			\psi(\bm{u}_t)-\psi(\bm{v}_{m_t})
		\right]
		=<\phi_i(\bm{r})\psi(\bm{v}),{\cal S}_3\{f^n\}-f^n>.
	\end{equation}
	Here we use Eq.~(\ref{eq_momentS3}). 
	On the other hand, the expected value of $Z_m$ is written as
	\begin{equation}
		E_p(Z_m)=\phi_i(\bm{r}_m)\left[\bar \psi-\psi(\bm{v}_{(m)})\right]
		\frac{\Delta t}{k}K[D_t\log S_{(m)}],
	\end{equation}
	and the expected value of $Z$ is written by the moment Eq.~(\ref{eq_momentfB})
	as
	\begin{equation}\label{eq_expZ}
		E_p(Z)=<\phi_i(\bm{x})\psi(\bm{v}),f^B>+{\cal O}(\Delta t^2).
	\end{equation}
	From Eqs.~(\ref{eq_realZ}) and (\ref{eq_expZ}), it is seen that the microscopic population density obtained by the tumbling process on $f^n$, ${\cal S}_3\{f^n\}$ is approximated by the sum of $f^n$ and $f^B$ under the assumption of the law of large numbers, i.e., $\hat Z\simeq E_p(Z)$,
	\begin{equation}\label{transS3}
		{\cal S}_3\{f^n\}=f^n+f^B+{\cal O}(\Delta t^2).
	\end{equation}

	In the division/death process ${\cal S}_4$, the particles in the $i$th lattice site may divide (or die) if the local macroscopic population density $\rho_i$ is smaller (or larger) than unity. 
	In the divisions, where $P[\rho_i]>0$, the particle creates a new particle with the same velocity at a random position within the same lattice site.
	Thus, when the divisions occur for $M_c$ particles for the distribution Eq.~(\ref{eq_distf}), the microscopic population density in the $i$th lattice site changes as
	\begin{equation}\label{eq_distS4_0}
		\begin{split}
			&{\cal S}_4\{f^n\}=\\
			&\frac{4\pi(\Delta x)^3}{M}\left\{
				\sum_{m=1}^{M_i^n}
				\delta(\bm{x}-\bm{r}^n_{(m)})
				\delta(\bm{v}-\bm{v}^n_{(m)})
				+\sum_{c=1}^{M_c}
				\delta(\bm{x}-(\bm{x}_{i+\frac{1}{2}}+\Delta x\bm{w}_c))
				\delta(\bm{v}-\bm{v}^n_{(m_c)})
			\right\},
		\end{split}
	\end{equation}
	where $\bm{w}$ is a random vector whose components are uniform random numbers in [$-\frac{1}{2}$,$\frac{1}{2}$].
	The moment is written as
	\begin{equation}\label{eq_momentS4_0}
		\begin{split}
			&<\phi_i(\bm{x})\psi(\bm{v}),{\cal S}_4\{f^n\}>\\
			&=\frac{1}{M}\left\{
				\sum_{m=1}^{M_i^n}\phi_i(\bm{r}_{(m)}^n)\psi(\bm{v}_{(m)}^n)
				+\sum_{c=1}^{M_c}\phi_i(\bm{x}_{i+\frac{1}{2}}+\Delta x\bm{w}_c)\psi(\bm{v}_{(m_c)}^n)\right\}.
			\end{split}
		\end{equation}
		In deaths, where $P[\rho_i]<0$, the particle is just removed from the lattice site.
		Thus, when the deaths occurs for $M_c$ particles in the $i$th lattice site, the microscopic population density obtained after the death process and its moment are written as changing the sign of the second term and replacing $\bm{x}_{i+\frac{1}{2}}+\Delta x\bm{w_c}$ with $\bm{r}_{(m_c)}^n$ in Eqs.~(\ref{eq_distS4_0}) and (\ref{eq_momentS4_0}), respectively.
		Hence, under the assumption of the uniform distribution of large number of particles in each lattice site, the moment for the population density after the division/death process ${\cal S}_4$ can be written as
		\begin{equation}\label{eq_momentS4}
			\begin{split}
				&<\phi_i(\bm{x})\psi(\bm{v}),{\cal S}_4\{f^n\}>\\
				&=\frac{1}{M}\left\{
					\sum_{m=1}^{M_i^n}\phi_i(\bm{r}_{(m)}^n)\psi(\bm{v}_{(m)}^n)
					+\mathrm{sign}(P[\rho_i])\sum_{c=1}^{M_c}\bar \phi_i\psi(\bm{v}_{(m_c)}^n)\right\} + {\cal O}(\frac{M_c}{M}\Delta x^2).
				\end{split}
			\end{equation}
			We remark that $M_c/M$ is estimated as $M_c/M\sim{\cal O}(\Delta t)$. 

			We introduce the stochastic variable $\Xi_m$ and $\Xi$ which are defined as
			\begin{equation}\label{eq_Xim}
				\Xi_m=\left\{
					\begin{array}{cl}
						\mathrm{sign}(P[\rho_i])\bar \phi_i \psi(\bm{v}_{(m)}),
						& (\mbox{if division/death}),\\
						0 &(\mbox{otherwise}),
					\end{array}
					\right.
				\end{equation}
				and
				\begin{equation}\label{eq_Xi}
					\Xi=\frac{1}{M}\sum_{m=1}^{M_i^n}\Xi_m.
				\end{equation}
				If the divisions(or deaths) occur for $M_c$ particles, the realized value of $\Xi$ is written as
				\begin{equation}
					\hat \Xi=\frac{\mbox{sign}(P[\rho_i])}{M}\sum_{c=1}^{M_c}
					\bar \phi_i
					\psi(\bm{v}_{(m_c)})
					=<\phi_i(\bm{x})\psi(\bm{v}),{\cal S}_4\{f^n\}-f^n>
					+{\cal O}(\Delta t\Delta x^2).
				\end{equation} 
				On the other hand, the expected value of $\Xi$ is written by Eq.~(\ref{eq_momentfC}) as
				\begin{equation}
					E_p(\Xi)=\frac{\Delta t P[\rho_i]}{M}\sum_{m=1}^{M_i^n}
					\bar \phi_i \psi(\bm{v}_{(m)})
					=<\phi_i(\bm{x})\psi(\bm{v}),f^C>+{\cal O}(\Delta t\Delta x^2).
				\end{equation}
				Thus, the microscopic population density obtained by the division/death process on $f^n$, ${\cal S}_4\{f^n\}$ is approximated by the sum of $f^n$ and $f^C$ under the assumption of the law of large numbers, i.e., $\hat \Xi\simeq E_p(\Xi)$,
				\begin{equation}\label{transS4}
					{\cal S}_4\{f^n\}=f^n+f^C+{\cal O}(\Delta t\Delta x^2).
				\end{equation}
				We remark that the error of the above equation is estimated at most ${\cal O}(\Delta t^2)$ when $\Delta t$ is at most the second order of $\Delta x$, i.e., $\Delta t\propto \Delta x^\alpha$ s.t. $\alpha\le 2$.

				In the MC simulation, the processes ${\cal S}_1$, ${\cal S}_3$, and ${\cal S}_4$ are successively conducted.
				For example, the process ${\cal S}_2$ is performed on the distribution obtained by the process ${\cal S}_1$ on $f^n$;
				\begin{equation}
					{\cal S}_3\{ {\cal S}_1\{f^n\} \}={\cal S}_3\{f^A\}=f^A+f^{B'},
				\end{equation}
				where $f^{B'}$ is obtained by replacing $f^n$ with $f^A$ in Eq.~(\ref{eq_fB}).
				However, the difference of $f^A$ and $f^n$ is ${\cal O}(\Delta t)$, so that $f^{B'}$ can be replaced with $f^B$ within the difference of ${\cal O}(\Delta t^2)$.
				Similarly, it is seen that the microscopic population density obtained by the successive three processes of ${\cal S}_1$, ${\cal S}_3$, and ${\cal S}_4$ approximates the microscopic population density at the next time step which satisfies the weak formulation Eq.~(\ref{eq_kinetic_dt}) within the difference of ${\cal O}(\Delta t^2)$, i.e.,
				\begin{equation}
					{\cal S}_4\{ {\cal S}_3\{ {\cal S}_1\{f^n\}   \} \}
					=f^A+f^B+f^C+{\cal O}(\Delta t^2)
					=f^{n+1}+{\cal O}(\Delta t^2).
				\end{equation}


\section{Concluding remarks}
We studied the self-organized pattern formation of chemotactic bacteria based on a kinetic chemotaxis model which includes a recently advocated formalism for bacterial chemotaxis, i.e., the logarithmic sensing of chemical cues along the pathway of bacterium and stiff and bounded signal response.
We have discovered a novel linear instability condition Eq.~(\ref{kinetic_instability}) stemming from the stiffness of chemotactic response.
Apart from the macroscopic description, we have been able to uncover the instability mechanism at the microscopic level.
The stationary homogeneous state of the macroscopic population density becomes linearly unstable and stationary periodic patterns are generated under the linear instability condition.
A remarkable property is that no eigenmodes exist in the large-oscillation limit in the linearized kinetic equation, which explains that pattern formations occur as observed in experiments.
Our new dispersion relation for instability also turns out to be sharp in the macroscopic limit, i.e., the flux-limited Keller-Segel equation.

MC simulations rigorously based on the kinetic chemotaxis model are performed with changing the parameters involved in the linear instability condition.
The numerical results demonstrate that the obtained linear instability condition is compatible and even sharply predicts the occurrence of the periodic pattern formations. See Fig.~\ref{fig_2dmap}.
The power spectra of the macroscopic population density show the plateau regime at the large wave numbers, where the perturbations neither grow nor damp irrespective of the linear instability condition. See Fig.~\ref{fig_turing_pow}
This observation is compatible with the fact that no eigenmodes exist in the large oscillation limit in the linearized kinetic chemotaxis equation.
Unexpectedly, the instability pattern undergoes transitions from the periodic oscillation around the uniform state $\rho=1$ to the periodic localized spikes over the zero-density state $\rho=0$ as increasing the modulation amplitude in chemotactic response. See Fig.~\ref{fig_spik}.

				\begin{table}[tbp]
					\centering
					\begin{tabular}{cc}
						\hline\hline
						$\psi_0$&3.0 [s]\\
						$t_0$&4500 [s]\\
						$V_0$&25 [$\mu$m/s]\\
						$\tilde d$&8$\times 10^{-6}$ [cm$^2$/s]\\
						$\tilde \delta^{-1}$& 20 [s]\\
						$\chi$ &0.2\\
						\hline\hline
					\end{tabular}
					\caption{Experimental values for mean tumbling frequency $\psi_0$, doubling time $t_0$, running speed $V_0$, diffusion coefficient of chemoattractant $\tilde d$, and stiffness $\tilde \delta^{-1}$ and modulation $\chi$ in chemotactic response obtained in Ref.~\cite{art:11SCBPBS}.}\label{t_experiment}
				\end{table}
The obtained instability condition Eq.~(\ref{kinetic_instability}) includes three control parameters, i.e., $k$, $d/k$, and $F'[0]/k$, which are written in the dimensional form as
				\begin{gather}
					k=1/(\psi_0t_0),\quad d/k=\tilde d\psi_0/(V_0^2t_0b_0),\quad F'[0]/k=\chi\tilde \delta^{-1}\psi_0.
				\end{gather}
Here $\psi_0$ is the mean tumbling frequency, $V_0$ is the running speed of bacteria, $b_0$ is the degradation rate of chemoattractant, and $t_0$ is the characteristic time which corresponds to the doubling time in cell division for Eq.~(\ref{eq_logistic}). 
See also the paragraph above Eq.~(\ref{kinetic_instability}).
Here $\tilde d$ and $\tilde \delta^{-1}$ are the diffusion coefficient of chemoattractant and stiffness of response function, respectively, in the dimensional form.
The values of control parameters are estimated from experimental data in Ref.~\cite{art:11SCBPBS} (See Table \ref{t_experiment}) as $k=7.5\times 10^{-5}$, $d/k=3.8\times (t_0b_0)^{-1}$, and $F'[0]/k=12$.
It is difficult to measure the degradation rate of chemoattractant $b_0$ in experiments.
In some references \cite{art:16EGBAV,art:06SZLJL,art:10SCBBSP}, the value of $b_0$ is estimated as $b_0=4\times 10^{-3}\sim 5\times 10^{-2}$ from the comparison between experimental and numerical results, so that we may estimate $d/k$ as $d/k\lesssim 1$.
Thus, from our analysis, we can expect the stationary homogeneous state becomes destabilized and pattern formation occurs for chemotactic bacteria.
Furthermore, for example, in Ref.~\cite{art:91BB}, it is argued that the pattern formation is suppressed by reduction of chemotactic sensitivity.
This argument is also consistent with our instability condition.
Although it remains to be assessed how quantitatively our instability condition explains the experimental results in terms of pattern formation, the present study convinces us that the self-organized pattern formation occurs due to the modulation of stiff response in chemotaxis in a realistic range of parameters.

Finally, our powerful MC method derived rigorously here has a possible advantage that can be extended to include internal states stemming from an intra-cellular chemical pathway. 
The kinetic chemotaxis model used in this study is based on a simplified model where the intra-cellular adaptation dynamics in chemotactic response is ignored and replaced by a instantaneous material derivative of chemical cue along the pathway of bacterium.
The time scale of adaptation, say $\tau_M$, is larger than the inverse of tumbling rate $\psi_0^{-1}$ but is much smaller than the doubling time $t_0$, i.e., $\psi_0^{-1}/t_0 < \tau_M/t_0 \ll 1$ \cite{art:05CG,art:12ZSDOWHJLT}.
Thus, the adaptation dynamics may be significant for the pattern formation.
In order to consider the adaptation dynamics internal states have to be taken into account in the tumbling kernel\cite{art:04EO,art:15X,art:16PTV}.
In a forthcoming paper, we plan to extend the MC method toward this direction and thus be able to challenge problems with another time scale related to the internal-state dynamics.


\section*{Acknowledgements} 
This study was financially supported by JSPS KAKENHI Grant Number 15KT0110 and 16K17554 and Institut Henri Poincar\'e RIP program.

\end{document}